\documentclass[conference]{IEEEtran}

\usepackage{graphics}
\usepackage{epsfig}
\usepackage{amsmath}
\usepackage{amssymb}
\usepackage[latin1]{inputenc}
\usepackage{amsfonts}
\usepackage{subfigure}
\usepackage{color}%
\usepackage{graphicx}
\providecommand{\U}[1]{\protect\rule{.1in}{.1in}}

\IEEEoverridecommandlockouts
\newtheorem{lemme}{Lemma}
\newtheorem{prop}{Proposition}

\title{{\LARGE \textbf{Fractional order differentiation by integration with Jacobi
polynomials}}}
\author{Da-Yan Liu, Olivier Gibaru, Wilfrid Perruquetti  and Taous-Meriem Laleg-Kirati \thanks{D.Y. Liu and T.M. Laleg-Kirati are with
Computer, Electrical and Mathematical Sciences and Engineering Division, King Abdullah
university of science and technology, KSA
\texttt{{\small Dayan.Liu@kaust.edu.sa; taousmeriem.laleg@kaust.edu.sa}}} \thanks{O. Gibaru is with LSIS
(CNRS, UMR 7296), Arts et Métiers ParisTech, 8 Boulevard Louis
XIV, 59046 Lille Cedex, France \texttt{{\small olivier.gibaru@ensam.eu}}%
}\thanks{W. Perruquetti is with LAGIS (CNRS, UMR 8146), École Centrale de
Lille, BP 48, Cit\'e Scientifique, 59650 Villeneuve d'Ascq, France
\texttt{{\small wilfrid.perruquetti@inria.fr}}}\thanks{O. Gibaru and W.
Perruquetti are with L'Équipe Projet Non-A, INRIA Lille-Nord Europe, 40,
Avenue Halley, 59650 Villeneuve d'Ascq, France.}}

\begin{document}

\maketitle
\thispagestyle{empty}
\pagestyle{empty}

\begin{abstract}
The differentiation by integration method with Jacobi polynomials was
originally introduced by Mboup, Join and Fliess
\cite{Mboup2007, Mboup2009a}. This paper generalizes this
method from the integer order to the fractional order for estimating the
fractional order derivatives of noisy signals. The proposed fractional order
differentiator is deduced from the Jacobi orthogonal polynomial filter and the
Riemann-Liouville fractional order derivative definition. Exact and simple
formula for this differentiator is given where an integral formula involving
Jacobi polynomials and the noisy signal is used without complex mathematical
deduction. Hence, it can be used both for continuous-time and discrete-time
models. The comparison between our differentiator and the recently introduced
digital fractional order Savitzky-Golay differentiator is given in numerical
simulations so as to show its accuracy and robustness with respect to
corrupting noises.
\end{abstract}


\section{INTRODUCTION} \label{section1}

Fractional models arise in many practical situations (\cite{Oustaloup10,Oustaloup02} for example). Such fractional order systems may also be used for control purposes: CRONE control is known to have good robustness properties (see \cite{Oustaloup95, Oustaloup98, Oustaloup99, Oustaloup02}). In order to implement such controller one needs to have a good digital fractional order differentiator from noisy signals, which is the scope of this paper.

The fractional derivative has a long history and are now very useful  in
science, engineering and finance (see, e.g., \cite{Oldham,Fliess97,Fliess98}).
The fractional order differentiator is concerned with estimating the
fractional order derivatives of an unknown signal from its noisy observed
data. Because of its importance, various methods have been developed during
the last years. They are divided into two kinds: continuous-time model (see,
e.g., \cite{Podlubny,Charef,Oustaloup}) and discrete-time model (see, e.g.,
\cite{Tseng1,Samadi,Tseng2}). Nevertheless, the real case of signals with
noises was somewhat overlooked. In order to resolve this problem, an
optimization formulation using genetic algorithms was proposed in
\cite{Machado} to reduce the noise effect in the estimations of the fractional
order derivatives. But, the complex mathematical deduction restricts its
application. A novel
\textbf{D}igital \textbf{F}ractional \textbf{O}rder \textbf{S}avitzky-\textbf{G}olay \textbf{D}ifferentiator (DFOSGD)
was introduced in \cite{Chen}, which was deduced from the Riemann-Liouville
fractional order derivative definition \cite{Podlubny2} (p. 63) and the
Savitzky-Golay filter \cite{Savitzky,Schafer}.
Using the Savitzky-Golay filter guarantees its accuracy and
simplicity for estimating the fractional order derivatives of noisy signals.
However, let us recall that the Legendre polynomial filter \cite{Persson} is
more recommended than the Savitzky-Golay filter for reasons of simplicity and
speed. In particular, it has advantages of being suitable for  irregularly
spaced or missing data.

The method of \emph{differentiation by integration} uses an integral of an
unknown noisy signal to estimate the integer order derivatives of this signal.
For free of noise signals, this method was firstly studied by
Cioranescu \cite{Cioranescu} (1938) and became well known for the Lanczos
generalized derivative \cite{C. Lanczos} (p. 324) (1956). The Lanczos
generalized derivative proposed an integral of a noisy signal and the Legendre
orthogonal polynomial to estimate the first order derivative of the signal.
Recently, in the noisy case, Mboup, Join and Fliess \cite{Mboup2007,
Mboup2009a} applied an algebraic setting to estimate high order derivatives by
integration, where Jacobi orthogonal polynomials were introduced in the
integral. Hence, we call the obtained differentiator \emph{Jacobi
differentiator}. Moreover, it was shown in \cite{Liu2011a} that the Jacobi
differentiator greatly improved the convergence rate of the Lanczos
generalized derivative. Very recently, it was shown in \cite{Liu2011c} that
the Jacobi differentiator could also be obtained by taking the derivative of
the Jacobi orthogonal polynomial filter considered as the extension of the
Legendre polynomial filter \cite{Persson,Meer}.

Let us recall that the algebraic differentiation method used in
\cite{Mboup2007, Mboup2009a} was introduced in \cite{FLI04a_compression_cras}
and also analyzed in \cite{Liu2009, Liu2011a, Liu2011b, Liu2011c, Riachy},
where the used algebraic manipulations were inspired by the algebraic
parametric estimation methods \cite{mexico, Mboup2009b, Liu2008, Liu2011d}.
Additional theoretical foundations can be found in \cite{ans,
FLI03_Indent_ESAIM}. In particular, by using non standard analysis techniques,
Fliess \cite{ans,shannon} showed that these methods exhibit  good robustness
properties with respect to corrupting noises without the need of knowing their
statistical properties. However, these methods have not been used to estimate
fractional order derivatives.

The aim of this paper is to generalize the differentiation by integration
method from the integer order to the fractional order for estimating
fractional order derivatives. In Section \ref{section2}, we recall the method
of integer order differentiation by integration with Jacobi polynomials. Then,
we deduce a fractional order differentiator from the Riemann-Liouville
fractional order derivative definition and the Jacobi orthogonal polynomial
filter. This differentiator is exactly given by an integral of Jacobi
polynomials. Hence, we call it \emph{fractional Jacobi differentiator}. In
Section \ref{section3}, we compare the fractional Jacobi differentiator to the
DFOSGD in some numerical simulations. Finally, we give some conclusions and
perspectives for our future work in Section \ref{section4}.


\section{METHODOLOGY}\label{section2}

Let $y^{\delta}=y+\varpi$ be a noisy signal observed in an
interval $I=[a,b] \subset\mathbb{R}$ of length $h=b-a$, where
$y \in\mathcal{C}(I)$ and the noise\footnote{More generally, the noise is a
stochastic process, which is bounded with certain probability and integrable
in the sense of convergence in mean square (see \cite{Liu2011b}).} $\varpi$ is
bounded and integrable. We are going to estimate the $\alpha^{th}$ ($\alpha
\in\mathbb{R}_{+} $) order derivative of $y$ by using its observation
$y^{\delta}$.


\subsection{Differentiation by integration with Jacobi polynomials}

In this subsection, we recall the method of integer order differentiation by
integration with Jacobi polynomials. This method was introduced in
\cite{Mboup2007,Mboup2009a} and studied in \cite{Liu2009, Liu2011a, Liu2011b,
Liu2011c, Riachy}.

The $n^{th}$ ($n\in\mathbb{N}$) order Jacobi orthogonal polynomial defined on
$[0,1]$ is given as follows (see \cite{Abramowitz} p. 775)
\begin{equation}
{P}_{n}^{(\mu,\kappa)}(t)=\displaystyle\sum_{j=0}^{n}\binom{n+\mu}{j}%
\binom{n+\kappa}{n-j}\left(  t-1\right)  ^{n-j}t^{j},\label{Eq_PolyJacobi}%
\end{equation}
where $\mu,\kappa\in]-1,+\infty\lbrack$. Let $f$ and $g$ be two functions
belonging to $\mathcal{C}([0,1])$, then we define the scalar product
$\left\langle \cdot,\cdot\right\rangle _{\mu,\kappa}$ of these functions by
(see \cite{Abramowitz} p. 774)
\begin{equation}
\left\langle f(\cdot),g(\cdot)\right\rangle _{\mu,\kappa}=\int_{0}^{1}{w}%
_{\mu,\kappa}(t)\,f(t)\,g(t)\,dt,
\end{equation}
where ${w}_{\mu,\kappa}(t)=(1-t)^{\mu}t^{\kappa}$ is the associated weighted
function. Thus, by denoting its associated norm by $\Vert\cdot\Vert
_{\mu,\kappa}$, we obtain
\begin{equation}%
\begin{split}
&  \Vert{P}_{n}^{\left(  \mu,\kappa\right)  }\Vert_{\mu,\kappa}^{2}=\\
&  \quad\quad\frac{\Gamma(\mu+n+1)\,\Gamma(\kappa+n+1)}{\Gamma(\mu
+\kappa+n+1)\,\Gamma(n+1)\,(2n+\mu+\kappa+1)}.
\end{split}
\label{Eq_Jacobinorm}%
\end{equation}
where $\Gamma(\cdot)$ is the classical Gamma function (see \cite{Abramowitz}
p. 255).

Let us ignore the noise for a moment. Then, we define an $N^{th}$ ($N
\in\mathbb{N}$) order approximation polynomial of $y$ by taking its truncated
Jacobi orthogonal series expansion
\begin{equation}
\label{Eq_Polyapprox1}%
\begin{split}
\forall\,t\in[0,1], \ \  &  D^{(0)}_{h,\mu,\kappa,N}y(a+h t):=\\
&  \sum_{i=0}^{N}\frac{\left\langle {P}_{i}^{(\mu,\kappa)}(\cdot
),y(a+h\cdot)\right\rangle _{\mu,\kappa}}{\Vert{P}_{i}^{(\mu,\kappa)}\Vert
^{2}_{\mu,\kappa}}\ {P}_{i}^{(\mu,\kappa)}(t).
\end{split}
\end{equation}
If we take $\kappa=\mu=0$, then the Jacobi orthogonal polynomials become the
Legendre orthogonal polynomials. Hence, this Jacobi polynomial filter
\cite{Meer} can be considered as the extension of the Legendre polynomial
filter \cite{Persson}.

If $y \in\mathcal{C}^n(I)$ with $n \in\mathbb{N}$, then we take the $n^{th}$ order derivative of the polynomial
$D^{(0)}_{h,\mu,\kappa,N}y(a+h \cdot)$ as an estimate of the $n^{th}$ order
derivative of $y$ \cite{Liu2011c}: $\forall\, t\in[0,1]$,
\begin{equation}
\begin{split}
& D_{h,\mu,\kappa,N}^{(n)}y(a+h t)\\ & \quad \quad := \frac{d^{n}}{d(ht)^{n}}\left\{ D_{h,\mu,\kappa,N}^{(0)} y(a+h
t)\right\} \\& \quad \quad =
\frac{1}{h^{n}} \frac{d^{n}}{dt^{n}}\left\{ D_{h,\mu,\kappa,N}^{(0)} y(a+h
t)\right\}.
\end{split}
\end{equation}
For any $t\in[0,1]$, this differentiator can be expressed as follows
\cite{Liu2011c}
\begin{equation} \label{Eq_Integer_Derivative_estimator1}
\begin{split}
&D_{h,\mu,\kappa,N}^{(n)}y(a+h t)=\\&\quad \quad \frac{1}{h^{n}} \int_{0}^{1} Q_{\mu
,\kappa,n,N}(\tau,t)\, y(a+h\tau)\, d\tau,
\end{split}
\end{equation}
with $C_{\mu,\kappa,n,i}=\frac{(\mu+\kappa
+2n+2i+1)\,\Gamma(\mu+\kappa+2n+i+1)\,\Gamma(n+i+1)}{\Gamma(\kappa
+n+i+1)\,\Gamma(\mu+n+i+1)}$,
$Q_{\mu,\kappa,n,N}(\tau,t) = w_{\mu,\kappa}(\tau)\displaystyle\sum
_{i=0}^{N-n} C_{\mu,\kappa,n,i} P_{i}^{(\mu+n,\kappa+n)}(t) {P}_{n+i}%
^{(\mu,\kappa)}(\tau)$.

Finally, we replace $y$ in (\ref{Eq_Integer_Derivative_estimator1}) by
$y^{\delta}$. Consequently, the $n^{th}$ order derivative of $y$ can be
estimated by an integral of Jacobi polynomials.
The noise effect on obtained estimations was analyzed in \cite{Liu2009, Liu2011b, Liu2011c, Mboup2009a, Mboup2009b}.
In the next subsection, we are going to generalize this differentiation by
integration method from the integer order to the fractional order.


\subsection{Fractional order differentiation by integration}

Similarly to \cite{Chen}, we are going to take the Riemann-Liouville
fractional order derivative of our approximation polynomial so as to get a
fractional order differentiator. The Riemann-Liouville fractional order
derivative (see \cite{Podlubny2} p. 62) is defined as follows: $\forall \, x\in\mathbb{R}_{+}^{\ast}$,
\begin{equation}
{}_{0}{D}_{x}^{\alpha}f(x):=\frac{1}{\Gamma(l-\alpha)}\frac{d^{l}}{dx^{l}}%
\int_{0}^{x}\left(  x-\tau\right)  ^{l-\alpha-1}%
f(\tau)\,d\tau,\label{Eq_RiemannLiouvillle}%
\end{equation}
where $0\leq l-1 \leq \alpha<l$ with $l\in\mathbb{N}^{\ast}$.

From now on, we denote the $\alpha^{th}$ ($\alpha\in\mathbb{R}_{+}
$) order derivative of $f$ by
\begin{equation}\label{}
\forall \, x\in\mathbb{R}_{+}^{\ast},\ {}_{0}{D}_{x}^{\alpha}f(x)=\frac{d^{\alpha}}{dx^{\alpha}}
\left\{f(x)\right\}= f^{(\alpha)}(x).
\end{equation}
Hence, if we take
$f(x)=x^{n}$ with $n\in\mathbb{N}$ and $x\in\mathbb{R}_{+}^{\ast}$, then we
obtain (see \cite{Podlubny2} p. 72)
\begin{equation}\label{Eq_Derivative_poly}%
\frac{d^{\alpha}}{dx^{\alpha}}\left\{x^{n}\right\}=\frac{\Gamma(n+1)}{\Gamma(n+1-\alpha
)}\,x^{n-\alpha}.
\end{equation}
By using (\ref{Eq_Derivative_poly}), we obtain the following lemma.

\begin{lemme}
\label{lemm_Derivative_Polyjacobi} The $\alpha^{th}$ ($\alpha\in\mathbb{R}_{+}
$) order derivative of the $n^{th}$ order Jacobi orthogonal polynomial
$P_{n}^{(\mu,\kappa)}$ defined in (\ref{Eq_PolyJacobi}) is given as follows:
$\forall \, t \in]0,1]$,
\begin{align}%
\begin{split}
& \frac{d^{\alpha}}{d t^{\alpha}}\left\{ P_{n}^{(\mu,\kappa)}(t)\right\}  =\\
& \quad\quad\quad\displaystyle\sum_{j=0}^{n}\sum_{l=0}^{n-j} c_{\mu
,\kappa,n,j,l}\,\frac{\Gamma(n-l+1)}{\Gamma(n-l+1-\alpha)} t^{n-l-\alpha},
\end{split}
\end{align}
where $c_{\mu,\kappa,n,j,l}= (-1)^{l}\binom{n+\mu} {j}\binom{n+\kappa}%
{n-j}\binom{n-j} {l}$.\newline
\end{lemme}

\noindent\textbf{Proof.} By applying the binomial theorem to
(\ref{Eq_PolyJacobi}), we get
\[
P_{n}^{(\mu,\kappa)}(t)=\displaystyle\sum_{j=0}^{n}\sum_{l=0}^{n-j}%
c_{\mu,\kappa,n,j,l}\,t^{n-l},
\]
where $c_{\mu,\kappa,n,j,l}=(-1)^{l}\binom{n+\mu}{j}\binom{n+\kappa}%
{n-j}\binom{n-j}{l}$. Hence, this proof can be completed by using
(\ref{Eq_Derivative_poly}) and the linearity of the fractional order
differentiation (see \cite{Podlubny2} p. 91). \hfill$\Box$

Then, we give the following proposition.

\begin{prop}
Let $y^{\delta}=y+\varpi$ be a noisy observation of $y$ in an interval
$I=[a,b] \subset\mathbb{R}$ of length $h=b-a$, where $y
\in\mathcal{C}(I)$ and the noise $\varpi$ is bounded and integrable. If the $\alpha^{th}$ ($\alpha\in
\mathbb{R}_{+}$) order derivative of $y$ exists, then a
fractional order differentiator, called \emph{fractional Jacobi
differentiator}, is defined as follows: $\forall \, t \in]0,1]$
\begin{equation}
\label{Eq_Fractional_Derivative_estimator}%
\begin{split}
&  D^{(\alpha)}_{h,\mu,\kappa,N}y^{\delta}(a+h t):=\\
& \quad\quad\quad\frac{1}{h^{\alpha}} \int_{0}^{1} Q_{\mu,\kappa,\alpha
,N}(\tau,t)\, y^{\delta}(a+h\tau)\, d\tau,
\end{split}
\end{equation}
where $h \in\mathbb{R}^*_{+}$, $N \in\mathbb{N}$, $\mu,\kappa\in]-1,+\infty
\lbrack$,
\begin{equation}
Q_{\mu,\kappa,\alpha,N}(\tau,t)=\sum_{i=0}^{N}w_{\mu,\kappa}(\tau) \frac
{P_{i}^{(\mu,\kappa)} (\tau)}{\Vert{P}_{i}^{(\mu,\kappa)}\Vert^{2}_{\mu
,\kappa}}q_{\mu,\kappa,i}(t),
\end{equation}
\begin{align}%
\begin{split}
q_{\mu,\kappa,i}(t)= \displaystyle\sum_{j=0}^{i}\sum_{l=0}^{i-j} c_{\mu
,\kappa,i,j,l}\,\frac{\Gamma(i-l+1)}{\Gamma(i-l+1-\alpha)} t^{i-l-\alpha},
\end{split}
\end{align}
with $c_{\mu,\kappa,i,j,l}= (-1)^{l}\binom{i+\mu} {j}\binom{i+\kappa}%
{i-j}\binom{i-j} {l}$ and $\Vert{P}_{i}^{(\mu,\kappa)}\Vert^{2}_{\mu,\kappa}$
is given in (\ref{Eq_Jacobinorm}).\newline
\end{prop}

\noindent\textbf{Proof.} Let us take the $\alpha^{th}$ order derivative of the
polynomial $D_{h,\mu,\kappa,N}^{(0)}y(a+h\cdot)$ defined in
(\ref{Eq_Polyapprox1}). By using the scale change property of the fractional order
differentiation (see \cite{Oldham} p. 76) we obtain: $\forall\,t\in]0,1]$,
\begin{equation}
\begin{split}
& D_{h,\mu,\kappa,N}^{(\alpha)}y(a+ht)\\& \quad \quad :=
  \frac{d^{\alpha}%
}{d(ht)^{\alpha}}\left\{  D_{h,\mu,\kappa,N}^{(0)}y(a+ht)\right\}\\
 & \quad \quad = \frac{1}{h^{\alpha}}\frac{d^{\alpha}%
}{dt^{\alpha}}\left\{  D_{h,\mu,\kappa,N}^{(0)}y(a+ht)\right\}
.\label{Eq_Derivative_estimator1}
\end{split}
\end{equation}
The linearity of the fractional order differentiation  gives us that
\begin{equation}%
\begin{split}
&  D_{h,\mu,\kappa,N}^{(\alpha)}y(a+ht)=\\
& \frac{1}{h^{\alpha}}\sum_{i=0}^{N}\frac{\left\langle {P}_{i}%
^{(\mu,\kappa)}(\cdot),y(a+h\cdot)\right\rangle _{\mu,\kappa}}{\Vert{P}%
_{i}^{(\mu,\kappa)}\Vert_{\mu,\kappa}^{2}}\frac{d^{\alpha}}{dt^{\alpha}%
}\left\{  P_{i}^{(\mu,\kappa)}(t)\right\}  .
\end{split}
\label{Eq_Fractional_Derivative_estimator2}%
\end{equation}
Finally, this proof can be completed by using Lemma
\ref{lemm_Derivative_Polyjacobi} and substituting $y$ by $y^{\delta}$ in (\ref{Eq_Fractional_Derivative_estimator2}). \hfill$\Box$

Consequently, by calculating the integral of the noisy signal $y^{\delta}$ and
a sum of the Jacobi polynomials we can estimate the value of $y^{(\alpha)}$ at
each point $a+ht$ in the interval $]a,b]$ for each $t\in]0,1]$. This integral
corresponds to a convolution in the discrete case.

If we fix the value of $t$ to $1$, then the fractional Jacobi differentiator
$D_{h,\mu,\kappa,N}^{(\alpha)}y^{\delta}(a+ht)$ only estimates the value of
$y^{(\alpha)}$ at point $b$. Thus, if we increase the length of the interval $I=[a,b]\subset
\mathbb{R}$, then we can estimate the other values of
$y^{(\alpha)}$. Consequently, the fractional Jacobi
differentiator can also be considered as a pointwise causal differentiator
which is useful in many on-line applications. Moreover,
we can fix the
parameters $\kappa$, $\mu$ and $N$ such that the sum of the Jacobi polynomials
can be explicitly calculated  by off-line work. Indeed, the fractional Jacobi
differentiator only needs a simple convolution for on-line applications.
The computation time is significatively improved.


\section{SIMULATIONS RESULTS}\label{section3}

The comparisons between the DFOSGD and some existing fractional order
differentiators were shown in \cite{Chen}, where the DFOSGD was better than
the others. In this section, by taking one numerical example considered in
\cite{Chen} we compare the fractional Jacobi differentiator to the DFOSGD so
as to show its accuracy and robustness with respect to corrupting noises.

\subsection{Numerical case}

From now on, we assume that $y^{\delta}(x_{i})=y(x_{i})+c\varpi(x_{i})$ is a
sequence of uniformly sampled noisy data of $y$ with $x_{i}=a+T_{s}i$ and
$T_{s}=\frac{h}{M}$ for $i=0,\cdots,M$. The noise $c\varpi(x_{i})$ is
assumed to  a zero-mean white Gaussian $iid$ sequence, and $c \in\mathbb{R}_{+}$. In this
discrete case, we apply the trapezoidal numerical integration method in the
fractional Jacobi differentiator to approximate the integral. Then, let us
recall that the values of $y^{(\alpha)}(x_{i})$ for $i=0,\cdots,M$ are
estimated by the DFOSGD as follows
\begin{equation}
\label{Eq_DFOSGD}\widetilde{y^{(\alpha)}(x_{i})}=\frac{1}{(T_{s})^{\alpha}}
b_{i,N} \left( X_{\theta,N}^{T} X_{\theta,N}\right) ^{-1}X_{\theta,N}^{T}
Y^{\delta}_{\theta},
\end{equation}
where $N \in\mathbb{N}$, $\theta\in\mathbb{N}$,
\begin{align}
b_{i,N}= \left(
\begin{array}
[c]{c}%
\frac{1}{\Gamma(1-\alpha)}\,(i+1)^{-\alpha}\\
\frac{1}{\Gamma(2-\alpha)}\,(i+1)^{1-\alpha}\\
\frac{\Gamma(3)}{\Gamma(3-\alpha)} \,(i+1)^{2-\alpha}\\
\vdots\\
\frac{\Gamma(N+1)}{\Gamma(N+1-\alpha)}\,(i+1)^{N-\alpha}\\
\end{array}
\right) ^{T}, \ \ Y^{\delta}_{\theta}= \left(
\begin{array}
[c]{c}%
y^{\delta}(x_{0})\\
y^{\delta}(x_{\theta})\\
y^{\delta}(x_{2\theta})\\
\vdots\\
y^{\delta}(x_{M})\\
\end{array}
\right) ,
\end{align}
and
\begin{align}
X_{\theta,N}= \left(
\begin{array}
[c]{cccc}%
1 & 1^{1} & \cdots & 1^{N}\\
1 & (1+\theta)^{1} & \cdots & (1+\theta)^{N}\\
1 & (1+2\theta)^{1} & \cdots & (1+2\theta)^{N}\\
\vdots & \vdots & \vdots & \vdots\\
1 & (M+1)^{1} & \cdots & (M+1)^{N}\\
&  &  &
\end{array}
\right) .
\end{align}
Hence, the idea is to obtain an $N^{th}$ order approximation polynomial by
using Savitzky-Golay algorithm \cite{Savitzky,Schafer} from
a subsequence $Y^{\delta}_{\theta}$, and then to calculate the $\alpha^{th}$
order derivative of this polynomial at each point $x_{i}$.


\subsection{Example}

\begin{figure}[ptb]
\centering {\includegraphics[scale=0.6]{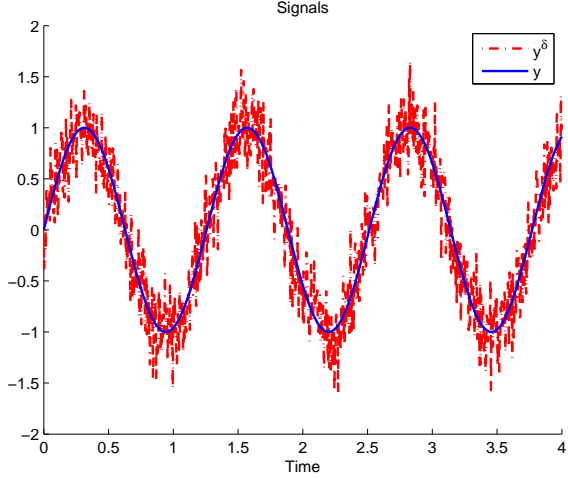}}\caption{ Signal $y$ and
noisy signal $y^{\delta}$.}%
\label{figure_signal}%
\end{figure}

In this example, we assume that $y^{\delta}(x_{i})=\sin(5x_{i})+c\varpi
(x_{i})$ where $x_{i} \in I= [0,4]$ for $i=0,\cdots,1000=M$, and $c=0.25$ is
adjusted in such a way that the signal-to-noise ratio $SNR=10\log_{10}\left(
\frac{\sum|y^{\delta}(x_i)|^{2}}{\sum|c\varpi(x_i)|^{2}}\right) $ is equal to
$SNR=10\text{dB}$. The noisy signal is shown in Figure \ref{figure_signal}.

\begin{figure}[ptb]
\centering {\includegraphics[scale=0.6]{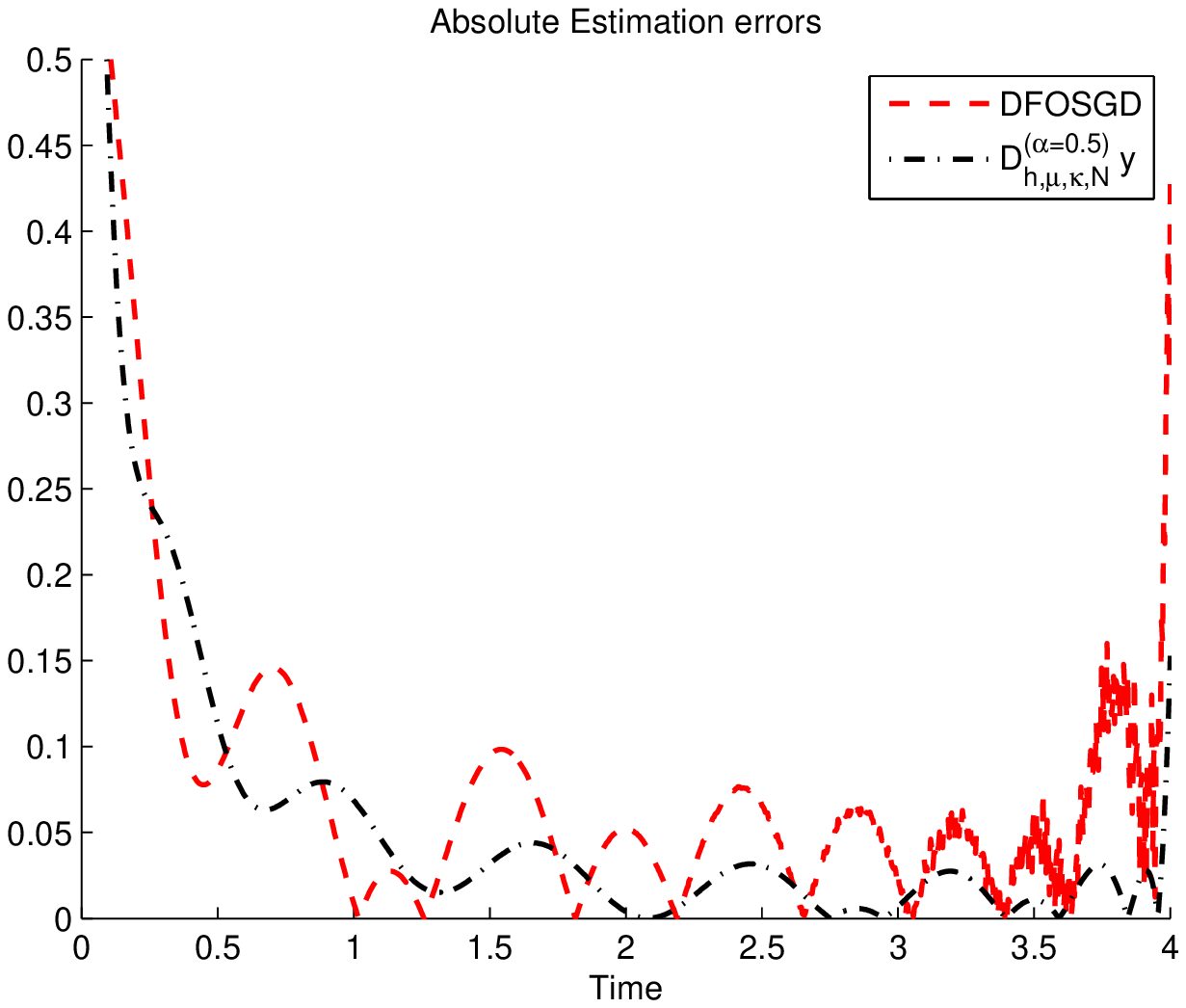}}\caption{Absolut
estimation errors in noise-free case.}%
\label{error1}%
\end{figure}

\begin{figure}[ptb]
\centering {\includegraphics[scale=0.6]{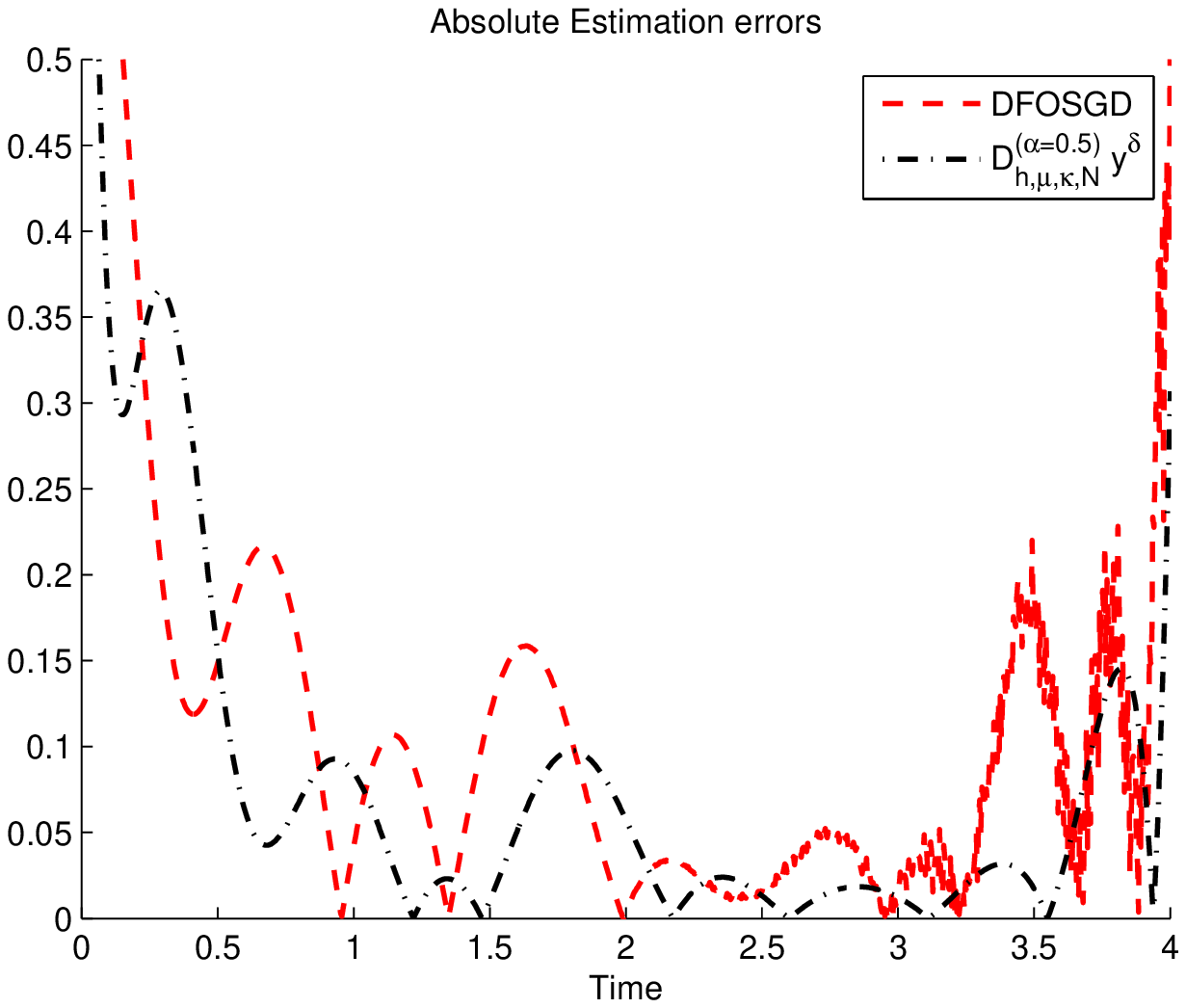}}\caption{Absolut
estimation errors in noisy case. }%
\label{error2}%
\end{figure}

\begin{figure}[ptb]
\centering {\includegraphics[scale=0.6]{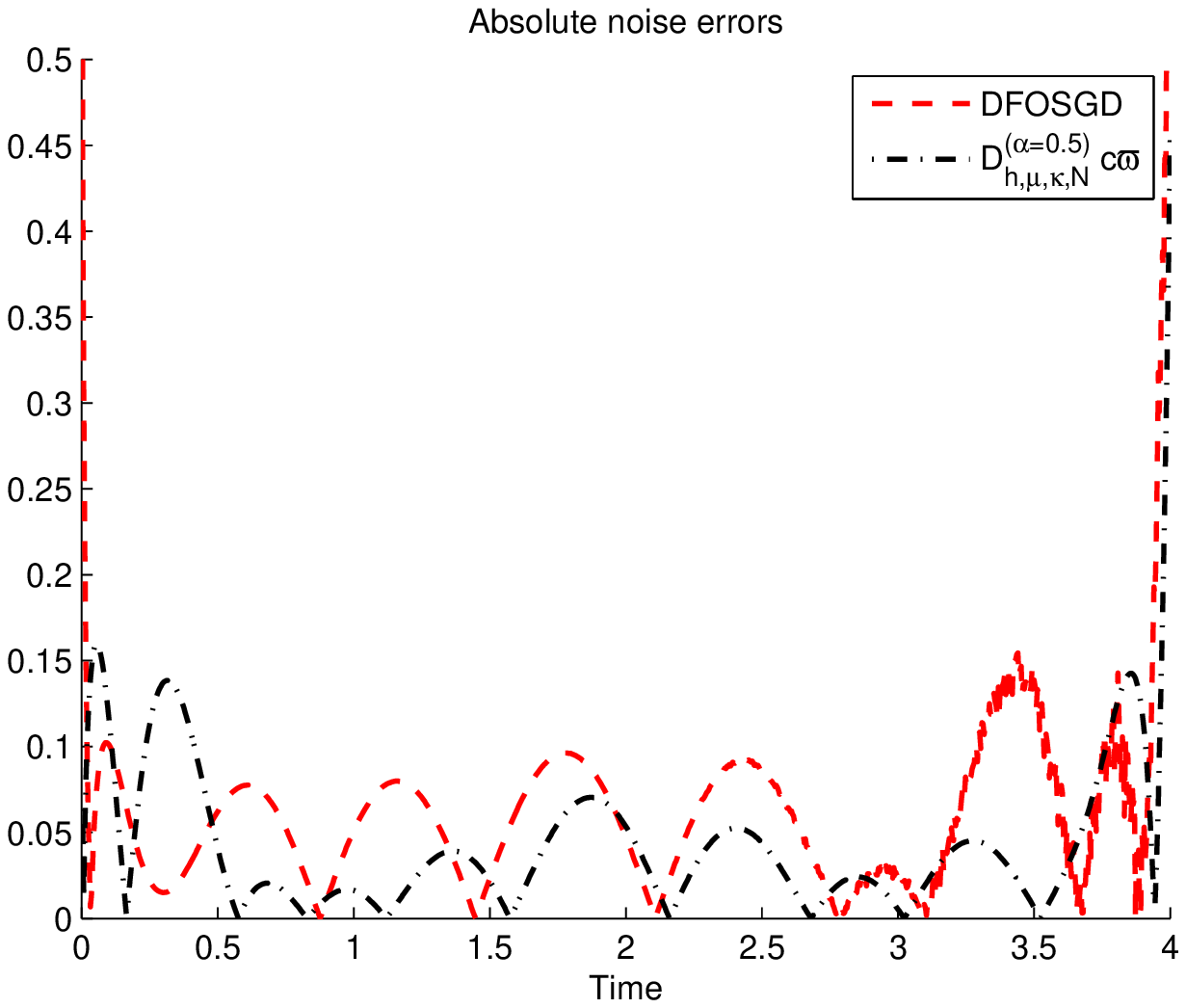}}\caption{Absolut noise
error contributions. }%
\label{noiseerror}%
\end{figure}

\begin{figure}[ptb]
\centering {\includegraphics[scale=0.6]{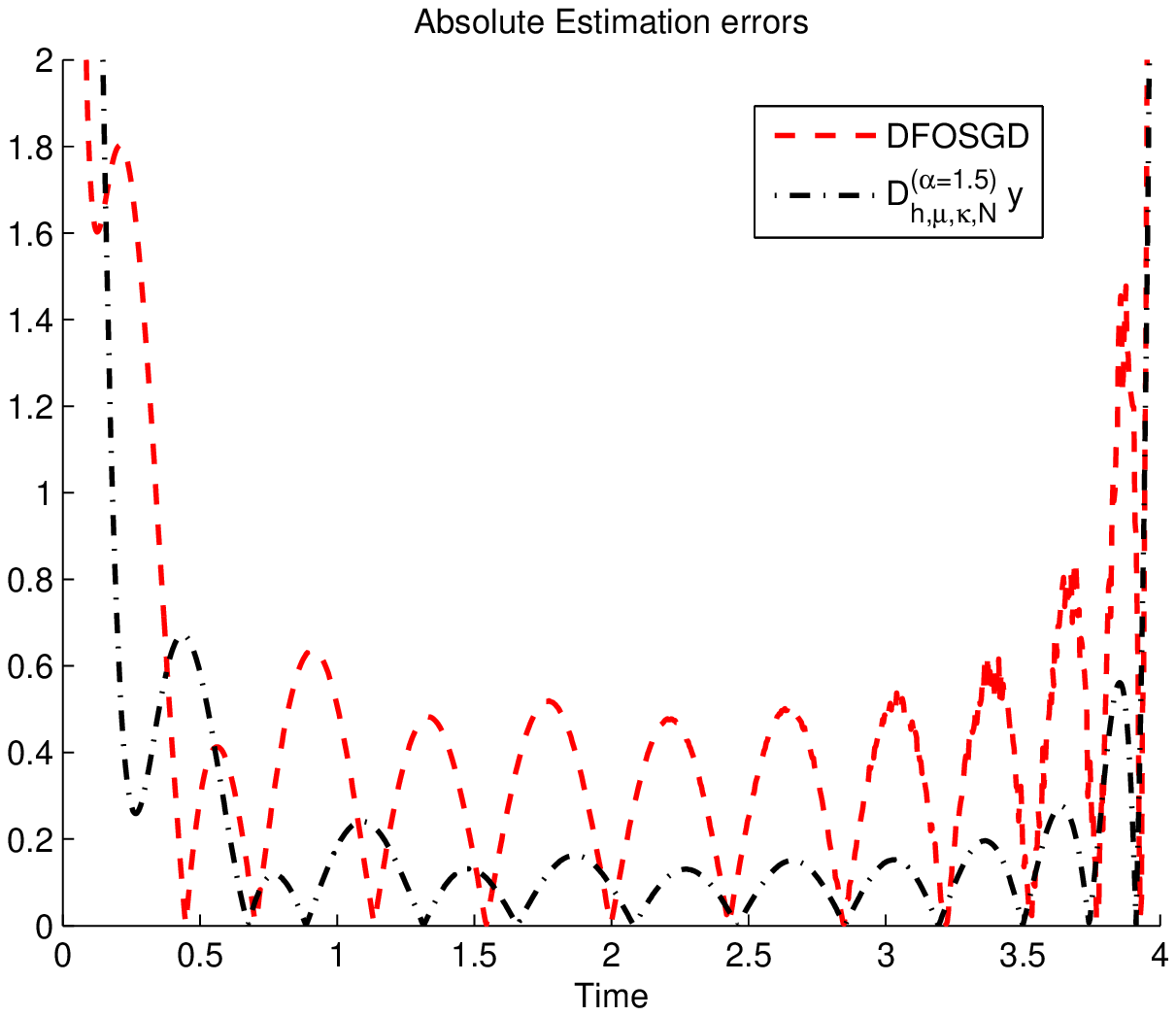}}\caption{Absolut
estimation errors in noise-free case.}%
\label{error3}%
\end{figure}

\begin{figure}[ptb]
\centering {\includegraphics[scale=0.6]{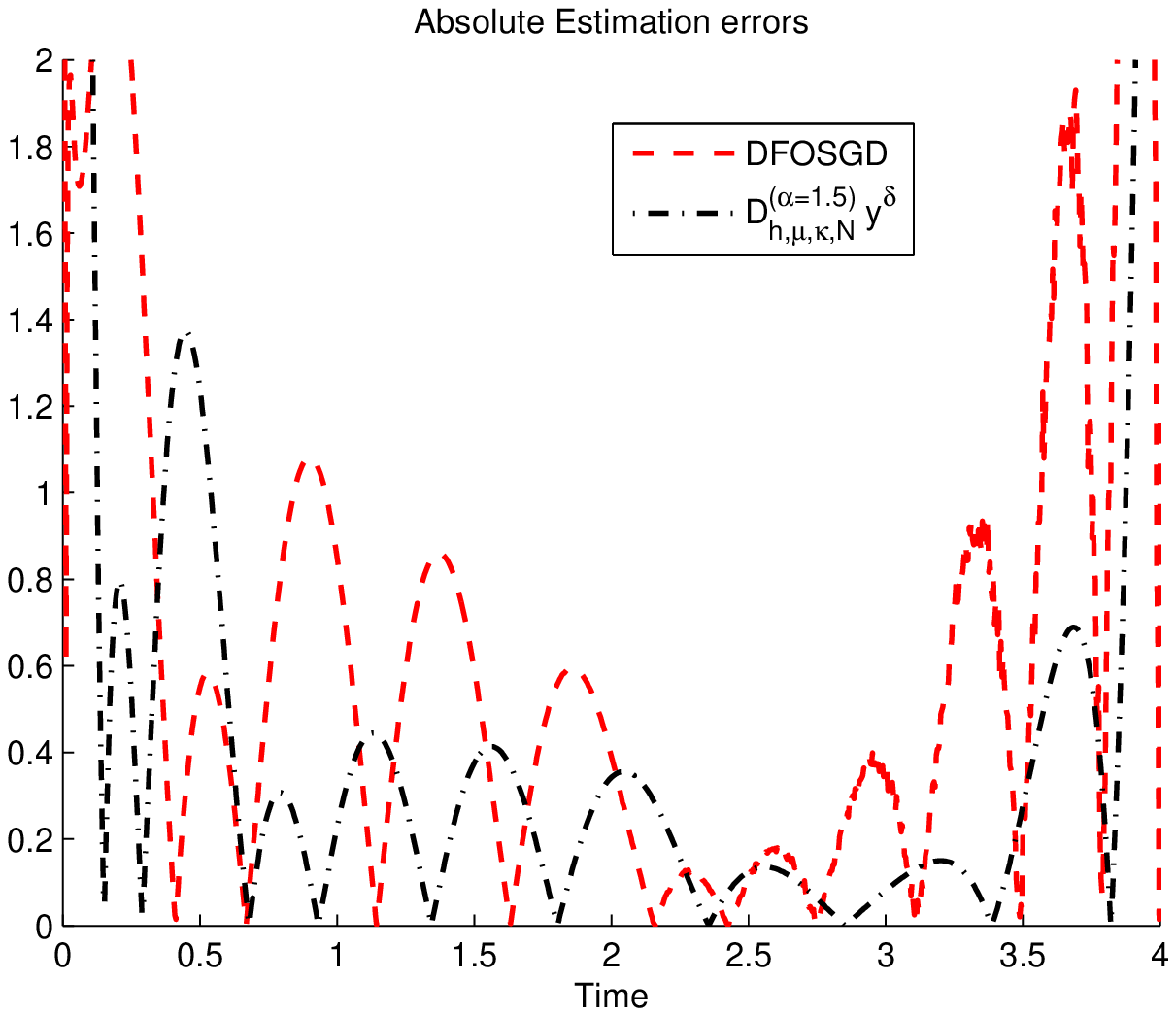}}\caption{Absolut
estimation errors in noisy case. }%
\label{error4}%
\end{figure}

\begin{figure}[ptb]
\centering {\includegraphics[scale=0.6]{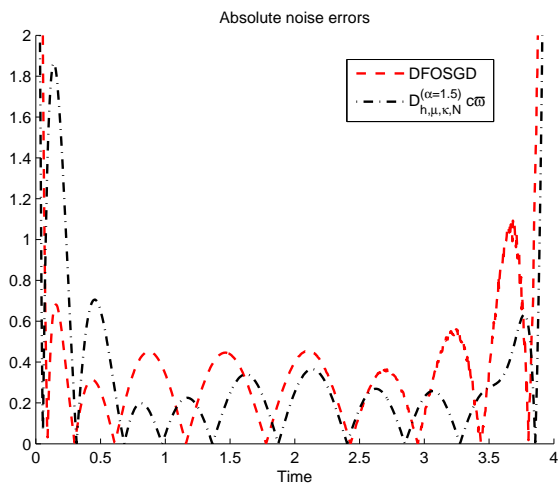}}\caption{Absolut noise
error contributions. }%
\label{noiseerror2}%
\end{figure}

Firstly, we  estimate the half order derivative of $y$. For the DFOSGD, we
take $N=14$ and $\theta=5$ which are the same values as
the ones used in \cite{Chen}. For the fractional Jacobi differentiator, we
take $N=14$ and $\kappa=\mu=0$. The associated absolute estimations errors in
the noise-free case and in the noisy case are given in Figure \ref{error1} and
Figure \ref{error2} respectively, where the black dash-dotted lines represent
the errors obtained by the fractional Jacobi differentiator and the red
dotted lines represent the ones obtained by the DFOSGD. Moreover, we can see
in Figure \ref{noiseerror} the associated absolute noise error contributions.
Secondly, we estimate the derivative $y^{(\alpha)}$ with $\alpha=1.5$ by using the
DFOSGD and the fractional Jacobi differentiator with $N=14$, $\theta=5$ and $\kappa=\mu=0$.
The obtained estimation errors are shown in Figure \ref{error3}, Figure \ref{error4} and Figure \ref{noiseerror2}.
Consequently, we can observe that the fractional Jacobi differentiator is more
accurate and more robust with noises than the DFOSGD in the noise-free case
and in the noisy case.

\section{CONCLUSION}\label{section4}

In this paper, we propose a fractional Jacobi differentiator
which is a fractional order differentiator exactly given by an integral
formula. This differentiator is deduced from the Jacobi orthogonal polynomial
filter and the Riemann-Liouville fractional order derivative definition. It
can accurately and easily estimate the fractional order derivatives of noisy
signals. There are some parameters $h$, $\mu$, $\kappa$ and $N$ on which the
fractional Jacobi differentiator depends. We do not give any analysis on the
influence of these parameters on the estimation errors. However, similarly to
the integer order differentiation by integration case (see
\cite{Liu2011b,Liu2011c,Liu2011a,Liu2009}), this study can be easily carried
out in a further work.  Moreover, we do not use the causal propriety of
the fractional Jacobi differentiator in the numerical simulations. It will be
used in some interesting applications where the sampling data may be
irregularly spaced.


\end{document}